\documentclass[3p,times]{elsarticle}

\usepackage{amssymb}
\usepackage{amsthm}
\usepackage{amsmath}

\usepackage{multirow}
 \biboptions{compress}
\usepackage{color}

\usepackage[figuresright]{rotating}
\usepackage{amssymb}
\newfont{\fp}{msbm10 at 11pt}

\usepackage[T1]{fontenc}
\usepackage{amsfonts}

\usepackage{amsfonts}
\usepackage{color, colortbl, framed}

\begin{document}

\begin{frontmatter}



\title{Structural break detection method based on the Adaptive Regression Splines
technique}


\author{Daniel Kucharczyk, Agnieszka Wy{\l}oma{\'n}ska}

\address{Hugo Steinhaus Center, Department of Mathematics,\\
    Wroclaw University of Technology \\
    Janiszewskiego 14a, 50-370 Wroc{\l}aw, Poland\\
    daniel.kucharczyk@pwr.edu.pl\\
    agnieszka.wylomanska@pwr.edu.pl}
\author{Rados{\l}aw Zimroz}
\address{ Diagnostics and Vibro-Acoustics Science Laboratory, \\
 Wroclaw University of Technology\\
 Na Grobli 15, 50-421 Wroclaw, Poland\\
		radoslaw.zimroz@pwr.edu.pl}
\begin{abstract}
For many real data, long term observation consists of different processes that coexist or occur one after the other.  Those processes very often exhibit  different statistical properties and thus before the further analysis the observed data should be segmented. This problem one can find in different applications and therefore new  segmentation techniques have been appeared  in the literature during last years. In this paper we propose a new method of time series  segmentation, i.e. extraction from the analysed vector of observations homogeneous parts with similar behaviour. This method is based on the absolute deviation about the median  of the signal and is an extension of the previously proposed techniques also based on the simple statistics. In this paper  we introduce the method of structural break point detection  which is based on  the Adaptive Regression Splines technique, one of the form of regression analysis. Moreover we propose also the  statistical test which allows testing hypothesis of behaviour related to different regimes. First, the methodology we apply to the simulated signals with different distributions in order to show the effectiveness of the new technique. Next, in the application part we analyse the real data set that represents the vibration signal from a heavy duty crusher used in a mineral processing plant.
\end{abstract}
\end{frontmatter}
\section{Introduction}
Physical variables observed/measured by advanced data acquisition systems in real world applications are often processed, modelled, and analysed in order to extract information about the process. In many cases, long term observation (long might mean very different period depends on application) in fact consists of different processes that coexist or occur one after the other. There are many examples that during single observation Process A reveals very different properties than Process B and they cannot be processed/analysed with the same tools (simple example is switching from stationary to non-stationary signal). Less radical case is that two segments might have the same type of model, but with different order and values of parameters. It leads to conclusion that very first step in experimental data analysis is signal segmentation, i.e. dividing the raw observation into smaller pieces (segments) with homogeneous properties. Note that Process A or Process B might be single process or mixture of processes.

\begin{figure}[h!]
\centering
\includegraphics[scale=0.35]{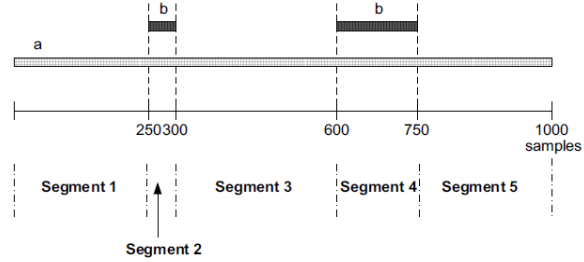}
\caption{The model of the signal that is a mixture of different segments related to different processes. Source: \cite{lopatka}.}
\label{lop}.
\end{figure}

From signal processing perspective, one of the most fundamental reason for segmentation is finding locally stationary segments in non-stationary signals. In such a context the segmentation is a problem of identifying the time instants at which the statistics of the observed signal changes. Signal segmentation is also generalized and considered in the context of signal detection and localization by means of recognition and identification of the time of appearance of event that differs from "normal" signal, \cite{lopatka}. It could be generalized to a model proposed in \cite{lopatka} presented in Fig. \ref{lop}.

The basis (criterion) for signal segmentation might be different in various domains, contexts and applications. Change of process from Process A to Process B might be related to appearance of seismic event, machine regime switching, fault in the system, financial crash on stock market, natural radiation, wind behavior, etc. A key issue is to find proper description of processes in order to highlight the difference between them. In some cases, the difference between processes is visually identifiable, but often advanced processing is required to find breakpoint.
Signal segmentation might be done in different manners. Some of the methods are based on the statistical properties of raw time series \cite{gsw,tsay} or	its modeling and analysis of model residuals \cite{lopatka,makowski1,makowski2}.
The second group of segmentation methods  takes under consideration not the raw signal but its  transformation to other domain, like spectral \cite{li}, cepstral \cite{kep}, time-frequency \cite{tf1,tf2,tf3}, etc. In the literature one can also find segmentation techniques based on the raw time series decomposition such as empirical mode or wavelets decomposition \cite{w1,w3,w2}.

Signal segmentation has been applied in many areas. It is especially crucial in condition monitoring (to isolate shocks related to damage)  \cite{tf2,w2}, machine performance analysis (to find when machine operates under overloading, idle mode etc) \cite{wyl_zimroz}, experimental physics \cite{gsw,bw,poiss}, biomedical signals (like ECG signals) \cite{kep,bio1,bio2,bio3,bio4,bio5,bio6}, speech analysis (automatic speech recognition and understanding) \cite{sp1,sp2,sp3}, econometrics \cite{eko1,eko2} and seismic signal segmentation \cite{tf1,sej1,sej2,sej3,sej4}.
The other areas where the segmentation problem appears one can find in \cite{ot1,ot2,ot3}.

In this paper we propose a novel  segmentation technique. This method is an extension of the algorithm proposed in \cite{gsw}, where the main statistics used to segmentation was based on the empirical second moment of raw data. An approach proposed here is a step forward in this field. The introduced method is based on the behaviour of the average absolute deviation about the empirical median of given time series. For the data where some statistical properties (expressed mostly by the scale parameter) change over time, the simple statistics applied here changes its behaviour which leads to structural break detection. The proposed method is much more sensitive for such cases where the change point is not visually identifiable in the contrast to the method presented in \cite{gsw}. Moreover it can be used to different raw signals without assumption of their distributions which is the main advantage. On the one hand  we analyse visually the mentioned statistics but on the other hand the strict method of change point detection is introduced. The method uses the Adaptive Regression Splines technique which is a form of regression analysis.  It is a non-parametric regression technique and can be seen as an extension of linear systems that automatically models non-linearities and interactions between variables \cite{spline}. Additionally, we propose a statistical test which allows for testing if given data satisfy the behaviour presented in Fig. \ref{lop}. Because in the analysis we consider the case where only one change point exists then  the time series $\{Z_i\}$ which is tested here, similar as in \cite{gsw}, can be written in the following mathematical form
\begin{equation}\label{eq:general}
    X_i = \left\{
            \begin{array}{ccc}
                T & \mbox{for} & i\leq l,\\\\
                Y & \mbox{for} & i> l
            \end{array}
            \right.
\end{equation}
where $l$ is a structural break point and $T$ and $Y$ are independent random variables with cumulative distribution functions $G(\cdot)$ and $ H(\cdot)$, respectively.

\begin{figure}[h!]
\centering
\includegraphics[scale=0.35]{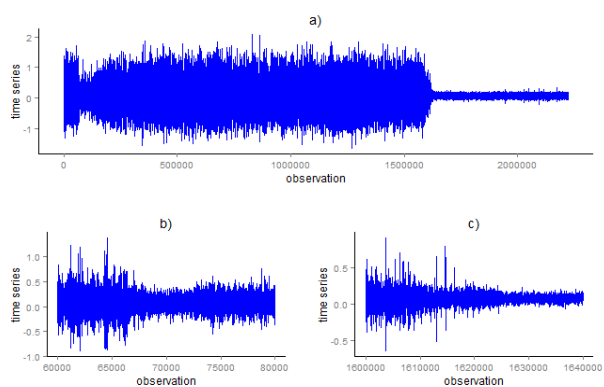}
\caption{The empirical time series  that represents the vibration signal from a heavy duty crusher used in a mineral processing plant (panel a)) and its subsignals corresponding to number of observations between $60000$ and $80000$ (panel b)) and $1600000$ and $1640000$ (panel c)).}
\label{fig:real_timeseries}.
\end{figure}

The motivation of our analysis is the time series presented in Fig. \ref{fig:real_timeseries} where we observe very specific behavior, similar to this described above. Namely, the visual inspection of the data indicates the statistical properties of raw signal change and there appears problem of structural break point detection, i.e. how to recognize automatically the changing point. The data presented in Fig. \ref{fig:real_timeseries}  represent the vibration signal from a heavy duty crusher used in a mineral processing plant. In the signal we observe that some statistical properties change over time which may be  related to the different modes of operation of the machine.  It was proved in many publications related to condition monitoring of machines, that vibration response significantly depends on external load applied to the machine.
In the case of crusher, load applied to machine means volume of material stream directed to crushing machine.
Depends on the volume as well as structure of material stream (size of stones) statistical properties of vibration will be different.
It is practically very hard to "quantify" properties of material stream.
By vibration signal segmentation we can find segments with homogeneous content that corresponds to "stable" or "stationary" material stream. Let us emphasize, observing the signal presented in Fig. \ref{fig:real_timeseries}, panel a) it seems that the structural break points are easy detectable here. However if we analyse separately the subsignals corresponding to number of observations between $60000$ and $80000$ (see panel b)) and $1600000$ and $1640000$ (see panel c))  we can conclude that in order to solve the  problem of structural break point detection we have to use advanced statistical methods.

The rest of the paper is organized as follows: In Section \ref{sec:methodology} we introduce the methodology related to the new segmentation technique and compare it to the method presented in \cite{gsw}. We propose here the visual test based on the  average absolute deviation about the empirical median of given time series as well as the new method of structural point detection based on the Adaptive Regression Splines technique. Next, we apply the proposed methodology to the various cases of simulated data with different distributions. We show the efficiency of the introduced method and their sensitivity for change point detection especially in case  where it is not visually  identifiable. We indicate the new technique can be also applied for data with heavy tailed distribution (like L\'evy stable) because it is based only on the empirical median of raw signal which always has finite value. In Section \ref{sec:RealDataAnalysis} we analyze the real data set presented in Fig. \ref{fig:real_timeseries} in the context of structural change point detection. At the end we recognize the distributions of segmented parts. Last section contains conclusions.
\section{Methodology}\label{sec:methodology}
\subsection{Regime testing- a scale based approach}\label{sec:scaleBased}
In this part, we introduce a new method of estimating the critical point of scale parameter for time series satisfying property given in \eqref{eq:general}. The general idea of the estimation procedure is closely linked to the approach introduced in \cite{gsw}, briefly described in the next part of this article, namely in Section~\ref{sec:GajdaEtAll}. In contrast to \cite{gsw}, our methodology is based on the behavior of the absolute deviation about median defined as follows
\begin{equation}\label{eq:NewStats}
V_j = \sum_{i=1}^{j} \left| X_i - \tilde{X}  \right|, \qquad j=1,2,\ldots,n,
\end{equation}
where $\tilde{X}$ is the median of time series $X_1, X_2, \ldots, X_n$. Under the assumption that the underlying sample $X_1, X_2, \ldots, X_n$ satisfies relation \eqref{eq:general}, we calculate the critical change point $l$ at which shift in the scale parameter occurs. It can be shown that in this case the estimated value of  the statistics \eqref{eq:NewStats} takes the following form
\begin{equation}\label{eq:basicMarsModel}
V_j = \beta_{0} + \beta_{1} \max{(0, j-l}) + \beta_{2}\max{(0, l-j)}.
\end{equation}
In the above formulation of the problem the procedure of estimating the critical change point $l$ is equivalent to the approach presented in \cite{gsw}.  One may easily proof  by performing simple algebra that equation (\ref{eq:basicMarsModel}) can be written as
\begin{equation}
V_j = \left\{\begin{array}{ccc}
a_1 + b_1j &\mbox{for}& j\leq l\\
a_2 + b_2j &\mbox{for}&  j> l\\
\end{array}\right. = \left\{\begin{array}{ccc}
\beta_0 -\beta_1l + \beta_1j &\mbox{for}& j\leq l\\
\beta_0 + \beta_2l -\beta_2j &\mbox{for}&  j> l\\
\end{array}\right.,
\end{equation}
where $a_1,b_1,a_2,b_2$ are the coefficients of the linear regressions fitted for two different segments of data separated by the point $l$ at which shift of the scale parameters occurs. \newline \indent For the purpose of estimating parameters $\beta_0, \beta_1, \beta_2$ defined in \eqref{eq:basicMarsModel} we use a non-parametric regression technique, called Multivariate Adaptive Regression Splines. For the sake of clarity we drop the description of the Multivariate Adaptive Regression Splines technique and describe it in the subsequent section, namely Section~\ref{sec:MARS}. For the purpose of testing of the occurrence of change defined in terms of the scale parameters of random variables $T$ and $Y$ in relation \eqref{eq:general}, we consider test proposed by \cite{AnsariBradley1960}, where the authors
test the hypothesis that $G(u)\equiv H(u)$ against alternatives of the form $G(u)\equiv F(\theta u)$ for $\theta\neq 1$.  Following the procedure described in \cite{AnsariBradley1960}, the two samples are ordered in a single joint array and ranks are assigned from each end of the joint array towards the middle. More formally, assuming that the two samples corresponding to random variables $T$ and $Y$ from relation \eqref{eq:general} consisting of $m$ and $n$ independent observations  are ranked in a combined array represented by
\begin{equation}\label{eq:ansariStats}
\mathcal{Z}=Z_1, Z_2, \ldots, Z_{m+n}.
\end{equation}
If $m+n$ then the vector of ranks $\mathcal{R}=\{\mathcal{R}_i\}$, $i=1,2,\ldots,m+n$ is as follows
\begin{equation}\label{eq:ranksEven}
\mathcal{R}=\{ 1,2,3\ldots,(m+n)/2, (m+n)/2\, \ldots, 3,2,1\}, \quad i=1,2,\ldots,m+n
\end{equation}
and, if $m+n$ is odd, the ranks are given in the following way
\begin{equation}\label{eq:ranksOdd}
\mathcal{R}=\{ 1,2,3\ldots,(m+n-1)/2, (m+n+1)/2, (m+n-1)/2,\ldots, 3,2,1\}, \quad i=1,2,\ldots,m+n.
\end{equation}
Finally, the test statistic $W$ under consideration is defined as the sum of the ranks given in \eqref{eq:ranksEven} or \eqref{eq:ranksOdd} associated with sample $T$. In the testing procedure we use the modified version of the statistics $W$, namely $W^{*}$ defined as
\begin{equation}\label{eq:ansariStatsW}
W^{*} = \left\{\begin{array}{ccc}
\displaystyle\frac{W-\frac{n(n+m+2)}{4}}{\sqrt{\displaystyle\frac{nm(n+m+2)(n+m-2)}{48(n+m-1)}}} & \mbox{ for } & $n+m=2k$\\
&&\\
\displaystyle\frac{W-\frac{n(n+m+1)^2}{4(n+m)}}{\sqrt{\displaystyle\frac{nm(n+m+1)(3+(n+m)^2)}{48(n+m)^2}}} & \mbox{ for } & $n+m=2k+1$\\
\end{array}\right., \quad k\in \mathbb{Z}.
\end{equation}
Under the null hypothesis, $G(u)\equiv H(u)$ the limiting distribution (where $m$ tends to $\infty$) of the statistics $W^{*}$ is characterized by a standard normal distribution $\mathcal{N}(0,1)$.
Therefore the corresponding $p$-value of statistics \eqref{eq:ansariStatsW} is calculated as $\Phi(W^*)$, where $\Phi(x)$ is defined as a cumulative distributive function of a standard normal distribution. If the $p$-value is equal to or smaller than the predefined significance level $\alpha$ we reject the null hypothesis $\mathcal{H}_0$, i.e. the time series under consideration does not satisfy the relation given in \eqref{eq:general}, in favor of the alternative hypothesis $\mathcal{H}_1$ stating that the observed time series has at least two regimes with different scale parameters.

\subsubsection{Multivariate Adaptive Regression Splines}\label{sec:MARS}
 Multivariate Adaptive Regression Splines (MARS) approach is a non-parametric, spline-based regression technique introduced in \cite{spline}. Let us consider a set of explanatory variables $X=(x_1,x_2,\ldots,x_n)$, a target variable $y$ and $N$ realizations of the process $\{y_i, x_{1i}, x_{2i}, \ldots x_{ni}\}_{1}^{N}$, whereas the true relationship between $X$ and $y$ is given in the form
\begin{equation}
\label{eq:trueRel}
y=f(x_1,x_2,\ldots, x_n) + \epsilon,
\end{equation}
where $f$ is an unknown function, and the error term $\epsilon$ is a white noise. The general idea of the MARS is to use the expansions in piecewise linear basis functions of the form $(x-t)_{+}$ and $(t-x)_{+}$, both called hockey stick basis functions with knot values set to $j$, where
\begin{equation}
\label{eq:hockeyStick}
(x-j)_{+} =
\left\{
\begin{array}{cc}
x-j,&\mbox{for } x>j\\
0,&\mbox{otherwise,}
\end{array}
\right.
\quad
\mbox{and}
\quad
(j-x)_{+} =
\left\{
\begin{array}{cc}
j-x,&\mbox{for } x<j\\
0,&\mbox{otherwise}
\end{array}
\right.
\end{equation}
\newline The general idea of the MARS is to form a set of pairs $\mathcal{C}$ consisting of basis functions given in \eqref{eq:hockeyStick} for each input variable $x_j$ with knots at each observed values $x_{ji}$. Therefore, the set of basis functions is given as
\begin{equation}\label{eq:setBasisFunctions}
\mathcal{C}=\left\{(x_j-t)_{+}, (t-x_j)_{+}\right\}\mbox{ for } t \in \{x_{j1},x_{j2},\ldots, x_{jN}\} \mbox{ and } j=1,\ldots n
\end{equation}
Under Multivariate Adaptive Regression Splines approach the set of basis functions defined in \eqref{eq:setBasisFunctions} is used to construct the final model as follows
\begin{equation}\label{eq:MARSModel}
f(X) = \beta_0 + \sum\limits_{m=1}^{M}\beta_mh_{m}(X),
\end{equation}
where each $h_m(X)$ is a function formed by taking single element from the initial set of basis functions $\mathcal{C}$, or a product of two or more such elements (functions). \newline\indent The way of estimating the parameters of the MARS model \eqref{eq:MARSModel} consists of two stages. At the first stage the MARS performs a forward stepwise regression procedure by considering at each stage a new basis pair from $\mathcal{C}$ including all possible interaction among explanatory variables and finally adding the pair of basis functions to the final model set $\mathcal{M}$ based on the residual squared error criterion. To make MARS approach computationally less expensive, the level of interactions among explanatory variables, as well as the maximum number of basis functions in model \eqref{eq:MARSModel} while performing the forward stepwise procedure are specified by the user. The forward stepwise regression procedure always starts by considering the initial model set $\mathcal{M}$ consisting only of the constant function $h_0(X)=1$. At the end of this procedure we have a large model of the form \eqref{eq:MARSModel} that typically overfits the data. The problem of overfitting the data is resolved by applying a backward stepwise regression procedure in the second step. In this step, the terms whose removal causes the smallest increase in residual squared error is deleted from the model $\mathcal{M}$ at each stage, producing an estimated best model $\hat{f}_{\lambda}$, where $\lambda$ denotes the number of terms defined in the considered model. The optimal value of $\lambda$ is calculated on the basis of the generalized cross-validation criterion defined as
\begin{equation}
\label{eq:MARSGCV}
GCV(\lambda)= \displaystyle\frac{\sum\limits_{i=1}^{N}(y_i-\hat{f}_\lambda(X))^2}{\left(1-M(\lambda)/N\right)^2},
\end{equation}
where $\left(1-M(\lambda)/N\right)^2$ is a penalty factor accounting for the increase variance resulting from a complex model and $M(\lambda)$ is the effective number of parameters, including both number of terms in the models and the number of parameters used for selecting the optimal location of the knots.
\subsection{Regime variance testing - a quantile approach}\label{sec:GajdaEtAll}
 In this subsection we briefly recall the general idea of the regime variance testing presented in \cite{gsw}. The procedure is based on the behaviour of the empirical second moment of a given one-dimensional sample $X_1, X_2, \ldots, X_n$. Namely, the authors consider statistics given in the following form
\begin{equation}\label{eq:GajdaEtAllStats}
C_j = \sum_{i=1}^{j} X_i^2,\qquad j=1,2,\ldots,n
\end{equation}
Under the assumption that the random variables $T$ and $Y$ defined in \eqref{eq:general} have distributions with finite second moments $\sigma_T^2$ and $\sigma_Y^2$, respectively, the statistics $C_j$ has the following property
\begin{equation}
\label{eq:EGajdaEtAllStats}
E(C_j) = \left\{
\begin{array}{ccc}
j\sigma_T^2 & \mbox{for} & j \leq l\\
\\
j\sigma_Y^2 & \mbox{for} & j > l.\\
\end{array}
\right.
\end{equation}
If $\sigma_T^2=\sigma_Y^2$, then the mean of $C_j$ statistics is equal to $j\sigma^2_T$ for all $j=1,2,\ldots, n$, therefore for independent identically distributed (i.i.d.) sample, the expected value of the statistics is a linear function with the shift parameter equal to zero. Similar behaviour of $E(C_j)$ is also observed under the milder assumption that the corresponding random variables $T$ and $Y$ given in \eqref{eq:general} have distributions with infinite second moments. \newline\indent
According to \cite{gsw}, the algorithm of estimating the critical point $l$ for the sample that fulfills relation \eqref{eq:general} starts with dividing for the fixed $k=1,2,\ldots,n$ the $C_j$ statistics into two disjoint sets $\{C_j:j=1,2,\ldots,k\}$ and $\{C_j:j=k+1,k+2,\ldots,n\}$. Next, the algorithm fits the linear regression lines for the first and the second set, respectively by minimizing the sum of squared error. Subsequently, the estimator $\hat{l}$ of the true change point $l$ is defined as the number $k$ that minimizes the following expression
\begin{equation}\label{eq:hatGajda}
\hat{l} = \arg\min\limits_{1\leq k \leq n}\left[\sum\limits_{j=1}^k\left(C_j - y_j^1(k)\right)^2 - \sum\limits_{j=k+1}^n\left(C_j - y_j^2(k)\right)^2\right]
\end{equation}
where $y_j^1(k)$ and $y_j^2(k)$ are the values of the linear functions fitted to the data belonging to the appropriate set of points $\{C_j:j=1,2,\ldots,k\}$ and $\{C_j:j=k+1,k+2,\ldots,n\}$.\newline
In the last step of the procedure, the underlying time series $X_1, X_2,\ldots, X_n$ is divided into two vectors $W_1 = [X_1^2,X_2^2, \ldots, X_l^2]$ and $W_2 = [X_{l+1}^2,X_{l+2}^2, \ldots, X_n^2]$, where quantiles ($q_{\alpha / 2}$ and $q_{1- \alpha/ 2}$) from the distributions of the squared time series from the vector $W_1$(for that the empirical standard deviation was smaller) are calculated according to the formula
\begin{equation}
P\left(q_{\alpha / 2} \leqslant X_i ^2 \leqslant 1-q_{\alpha / 2} \right) = 1-\alpha, \qquad i=1,2,\ldots, l
\end{equation}
where $\alpha$ is a given significance level. The above formula is used to construct the $\mathcal{H}_0$ hypothesis defined as previously: observed time series does not satisfy relation \eqref{eq:general}, whereas the alternative hypothesis $\mathcal{H}_1$ is formulated as:
observed time series has at least two regimes of the data for which
the appropriate quantiles of the squared time series are different.
Under the $\mathcal{H}_0$ hypothesis, the appropriate quantiles can be determined on the basis of the empirical cumulative distribution function. Due to the fact that $X_{l+1}^2,X_{l+2}^2, \ldots, X_n^2$ are independent, the statistics $B$ given as
\begin{equation}\label{eq:quantileTestStats}
B = \#\left\{q_{\alpha / 2} \leqslant W_2 \leqslant 1-q_{\alpha / 2}\right\}
\end{equation}
has a binomial distribution with parameters $n-l$ and $p=1-\alpha$.
\subsection{Methodology overview}
\begin{figure}[h!]
\centering
\includegraphics[scale=0.4]{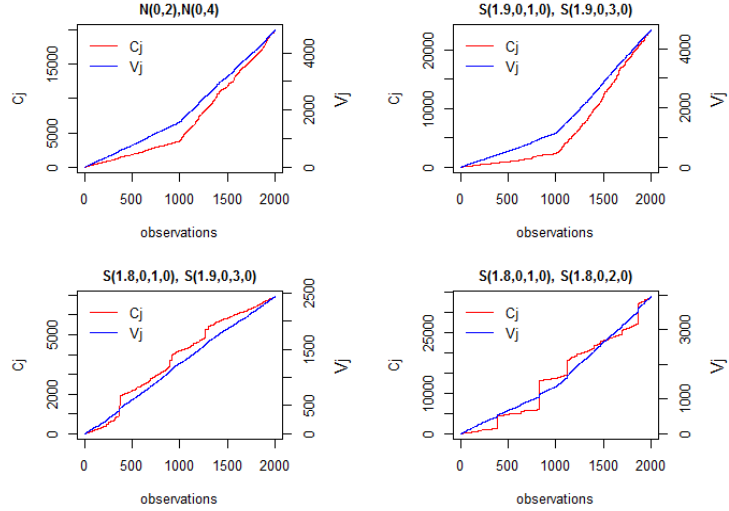}
\caption{The $C_j$ (red line) and $V_j$ (blue line) statistics defined in \eqref{eq:GajdaEtAllStats} and \eqref{eq:NewStats} for time series driven by a combination of distributions with different set of parameters where the corresponding critical point $l$ has been planned at middle of the length of sample.}
\label{fig:StatsGraph}
\end{figure}

 This section is supposed to describe the main differences between the methodologies described in the foregoing sections, see Sections~\ref{sec:scaleBased} and ~\ref{sec:GajdaEtAll}. The main difference comes from the problem of detecting an appropriate change point at which the shift in the dispersion of a distribution has taken place is tackled by using two different measures of variability of the distribution. Namely, the average absolute deviations about median and the second moment of a given sample have been considered. The former is used as a direct measure of scale parameters of a distribution with cumulative distribution function $F$ satisfying relation
\begin{equation}
F(x;s,\theta) = F(x/s;1,\theta),
\end{equation}
where $s$ is a scale parameter. The latter one is defined in a more systematic way by measuring the variability of distribution in terms of the empirical variance. As it appears that the graphs of diagnostic statistics presented in \eqref{eq:NewStats} and \eqref{eq:GajdaEtAllStats}  reveal the same "hockey stick" pattern for time series satisfying relation \eqref{eq:general}. However the shape of the function derived based on \eqref{eq:NewStats} may contain points at which sudden jump (discontinuity)  is observed, influencing the bias and variance of the true critical point estimator $\hat{l}$ defined in $\eqref{eq:hatGajda}$. From the other hand, the resulting graph of the statistics defined in  \eqref{eq:NewStats} is smoother than in case of diagnostic statistics considered in \eqref{eq:GajdaEtAllStats}. As a consequence the obtained results of estimating the true critical point are more reliable comparing to those results obtained based on the behavior of the empirical second moment taken into consideration in \cite{gsw}. More details about the performance of the aforementioned statistics can be found in the successive part of this article, see Section~\ref{sec:simulations}.
Moreover the behaviour of both considered statistics, namely $V_j$ and $C_j$ for exemplary distributions of random variables $Z$ and $Y$ in relation \eqref{eq:general} is presented in Fig.~\ref{fig:StatsGraph}.
\section{Simulations}\label{sec:simulations}
In this section we describe the results of performed simulations for the methods presented in Section \ref{sec:methodology} and compare the robustness of detecting the critical point of the underlying sample satisfying relation \eqref{eq:general}. The comparison is made by simulating 1000 trajectories of length $n=1800$ of stochastically independent random variables with the scale parameter change placed on 800th observation. For the purpose of simulations we consider two different set of scenarios of family of distributions as in \cite{gsw}. For the first scenario, where the difference between distributions of random variables $T$ and $Y$ in equation (\ref{eq:general}) is relatively small, we consider the following cases
\begin{itemize}
\item the pure Gaussian case with $\mathcal{N}(0,4)$ and $\mathcal{N}(0,4.55)$ distribution for the first 800 and last 1000 observations, respectively.
\item the pure L\'evy-stable case with $\mathcal{S}(1.9,0,2,0)$ and $\mathcal{S}(1.9,0,2.5,0)$ distribution for the first 800 and last 1000 observations, respectively.
\item the pure L\'evy-stable case with $\mathcal{S}(1.8,0,2,0)$ and $\mathcal{S}(1.85,0,2.5,0)$ distribution for the first 800 and last 1000 observations, respectively.
\item the pure L\'evy-stable-Gaussian case with $\mathcal{S}(1.8,0,1.2,0)$ and $\mathcal{N}(0,2.45)$ distribution for the first 800 and last 1000 observations, respectively.
\end{itemize}
In the second scenario, we consider following parameters of distributions
\begin{itemize}
\item the pure Gaussian case with $\mathcal{N}(0,2)$ and $\mathcal{N}(0,4)$ distribution for the first 800 and last 1000 observations, respectively.
\item the pure L\'evy-stable case with $\mathcal{S}(1.9,0,2,0)$ and $\mathcal{S}(1.9,0,4,0)$ distribution for the first 800 and last 1000 observations, respectively.
\item the pure L\'evy-stable case with $\mathcal{S}(1.85,0,2,0)$ and $\mathcal{S}(1.95,0,4,0)$ distribution for the first 800 and last 1000 observations, respectively.
\item the pure Gaussian-L\'evy-stable case with $\mathcal{N}(0,4)$ and $\mathcal{S}(1.9,0,1,0)$  distribution for the first 800 and last 1000 observations, respectively.
\end{itemize}
Comparing the set of scenarios taken into account in our simulation study with those considered in \cite{gsw}, two additional test cases have been added in order to test specific behaviour of time series extracted from the real data set presented in Section~\ref{sec:RealDataAnalysis}. Namely, the newly introduced test cases exhibit the behaviour where the shift in the stability parameter $\alpha$ in conjunction with the change of the scale parameter has taken place.

\begin{figure}[h!]
\centering
\includegraphics[scale=0.35]{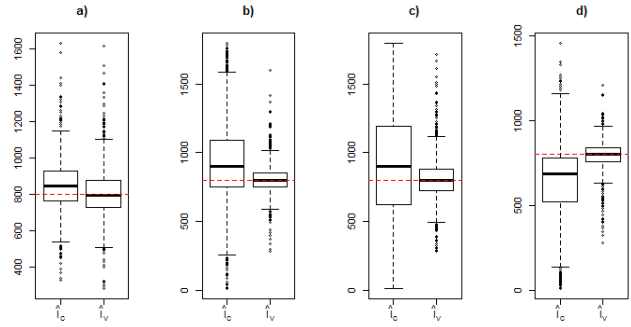}
\caption{Comparison of detection procedure for the critical scale change point for two estimators $\hat{l}_C$ and $\hat{l}_{V}$. Panel a) $\mathcal{N}(0,4),\  \mathcal{N}(0,4.55)$ ,panel b) $S(1.9,0,2,0),\ S(1.9,0,2.5,0)$, panel c) $S(1.8,0,2,0),\ S(1.85,0,2.5,0)$, panel d) $S(1.8,0,1.2,0),\  \mathcal{N}(0,2.45)$}
\label{fig:boxPlotNV}
\end{figure}

\begin{figure}[h!]
\centering
\includegraphics[scale=0.35]{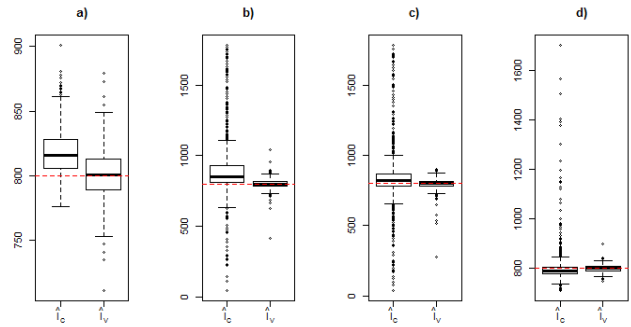}
\caption{Comparison of detection procedure for the critical scale change point for two estimators $\hat{l}$ and $\hat{l}_{new}$. Panel a) $\mathcal{N}(0,2),\  \mathcal{N}(0,4)$, panel b) $S(1.9,0,2,0),\ S(1.9,0,4,0)$, panel c) $S(1.85,0,2,0),\ S(1.95,0,4,0)$, panel d) $\mathcal{N}(0,4),\ S(1.9,0,1,0) $}
\label{fig:boxPlotV}
\end{figure}

\begin{figure}[h!]
\centering
\includegraphics[scale=0.3]{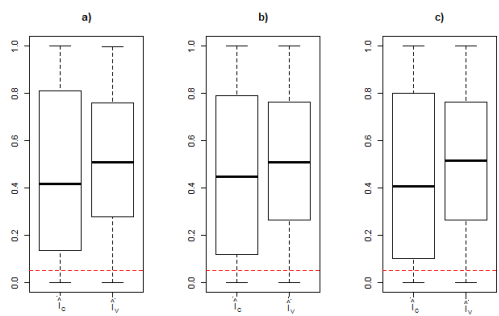}
\caption{The boxplots of p-values from selected test statistics obtained for different data sets, with the constant scale parameters. Panel a) $\mathcal{N}(0,2)$, panel b) $S(1.8,0,1.2,0)$, panel c) $\mbox{permuted } \mathcal{S}(1.8,0,1,0), \mathcal{S}(1.9,0,1,0)$. The corresponding significance level $\alpha=0.05$ is shown as red dashed line. }
\label{fig:TypeIErrorPValue}
\end{figure}

\noindent In terms of the results for the first scenario under consideration, where all test cases are based on the assumption that the corresponding critical point is hardly visible are presented in Fig.~\ref{fig:boxPlotNV}. The obtained results for the aforementioned scenario show that the estimator $\hat{l}_{V}$ of a true change point $l$ defined in \eqref{eq:MARSModel} outperforms the corresponding estimator $\hat{l}_{C}$ (calculated based on  \eqref{eq:hatGajda}) for all considered test cases. The results for the second scenario are presented in Fig.~\ref{fig:boxPlotV}, accordingly. Also in this case we can observe that $\hat{l}_{V}$ estimates the true change point in a more stable way than $\hat{l}_{C}$. The superiority of $\hat{l}_{V}$ is particularly pronounced in case of the pure L\'evy-stable case with  $\mathcal{S}(1.9,0,2,0)$ and $\mathcal{S}(1.9,0,4,0)$ distribution, where the dispersion of $\hat{l}_{V}$ is much more smaller than dispersion of $\hat{l}_{C}$.
Apart from comparing the detection procedures for the critical scale change point given above, we also examine the performance of the proposed estimator \eqref{eq:NewStats} in conjunction with test statistic by checking the type I error (rejection of  $\mathcal{H}_0$ hypothesis in case it is true) for 1000 Monte Carlo trajectories of length of 1800 of stochastically independent random variables for the following cases
\begin{itemize}
\item Gaussian case with $\mathcal{N}(0,2)$ distribution,
\item L\'evy-stable case with $\mathcal{S}(1.8,0,1.2,0)$ distribution
\item L\'evy-stable case with $\mathcal{S}(1.8,0,1,0)$ and $\mathcal{S}(1.9,0,1,0)$ distribution for each half of the sample, randomly permuted.
\end{itemize}
The results of the conducted simulations are presented in Table~\ref{tbl:TypeIerror}. For the testing purposes, we apply the sample mean values of obtained estimators $\hat{l}_{C}$ and $\hat{l}_{V}$ presented in Sections \ref{sec:GajdaEtAll} and \ref{sec:scaleBased}, respectively. Furthermore we assume that the corresponding significance level $\alpha$ is set to $0.05$.  The test based on the $C$ statistics  we denote as $GSW-C$ while the test proposed in this paper as $KWZ-V$.
\begin{table}[h!]
\centering
\begin{tabular}{c|c|c}
Distribution of sample & GSW-C &  KWZ-V \\\hline\hline
$\mathcal{N}(0,2)$ & 157 & 45  \\
$\mathcal{S}(1.8,0,1.2,0)$& 146 & 49 \\
$\mbox{permuted } \mathcal{S}(1.8,0,1,0), \mathcal{S}(1.9,0,1,0)$ & 168 & 51
\end{tabular}
\caption{Numbers of the incorrect rejection of a true null hypothesis $\mathcal{H}_0$ (type I error) based on 1000 Monte Carlo trajectories and the significance level $\alpha=0.05$.}
\label{tbl:TypeIerror}
\end{table}\newline
For the three considered cases, we see that the approach proposed in \cite{gsw} does not preserve the assumed type I error ($\alpha = 0.05$), as the corresponding type I errors are equal to $0.157$,  $0.146$ and $0.168$. In contrast to \cite{gsw}, the resulting type I errors for the test statistics described in Section \ref{sec:scaleBased} are preserved, see $KWZ-V$ column.
Moreover, the $p$-values corresponding to the acceptance of $\mathcal{H}_0$ are much more higher than significance level $\alpha=0.05$, see Fig.\ref{fig:TypeIErrorPValue}.
The next goal of our simulations study is to exploit the statistical power of the examined test and compare it with the reference test proposed in \cite{gsw}. As the power of a given test is strictly connected to type II error $\beta$ of the test, we focus our attention only on the values of the corresponding type II errors, i.e. numbers of wrongly accepted a false null hypothesis $\mathcal{H}_0$. For this reason, we consider the first set of scenarios of family of distributions described at the beginning of this section.
\begin{table}[h!]
\centering
\begin{tabular}{c|c|c}
Distribution of sample & GSW-C &  KWZ-V \\\hline\hline
$\mathcal{N}(0,4), \mathcal{N}(0,4.55)$ & 257 & 164  \\
$S(1.9,0,2,0),\ S(1.9,0,2.5,0)$& 219 & 0\\
$S(1.8,0,2,0),\ S(1.85,0,2.5,0)$& 540 & 1\\
$S(1.8,0,1.2,0),\ \mathcal{N}(0,2.45)$& 297 & 0
\end{tabular}
\caption{Numbers of the incorrectly accepted a false null hypothesis $\mathcal{H}_0$ (type II error) based on 1000 Monte Carlo trajectories and the significance level $\alpha=0.05$.}
\label{tbl:TypeIIerror}
\end{table}
\begin{figure}[h!]
\centering
\includegraphics[scale=0.35]{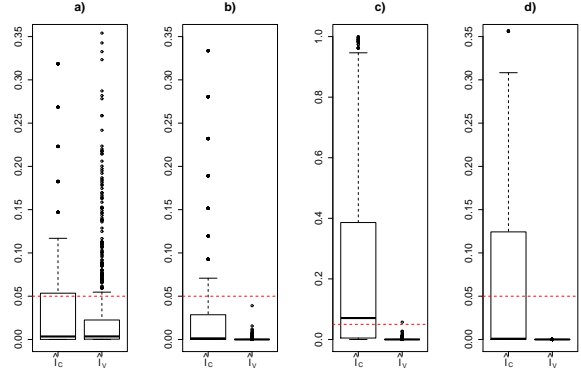}
\caption{The boxplots of p-values for different set of distributions, where initial change point is hardly visible. Panel a) $\mathcal{N}(0,4), \mathcal{N}(0,4.55)$, panel b) $S(1.9,0,2,0),\ S(1.9,0,2.5,0)$, panel c) $S(1.8,0,2,0),\ S(1.85,0,2.5,0)$, panel d) $S(1.8,0,1.2,0),\ \mathcal{N}(0,2.45)$. The corresponding significance level $\alpha=0.05$ is depicted as red dashed line.}
\label{fig:TypeIIErrorPValue}
\end{figure}
\newline
\newpage
\section{Real data analysis}\label{sec:RealDataAnalysis}
\begin{figure}[h!]
\centering
\includegraphics[scale=0.35]{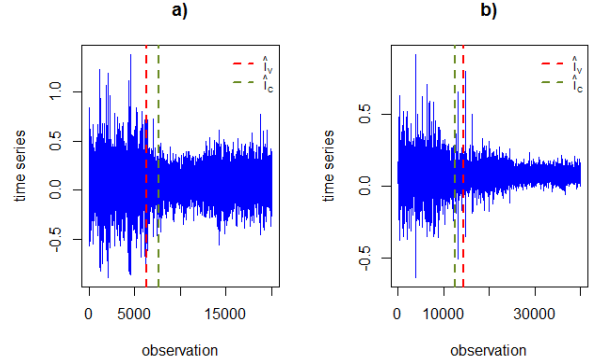}
\caption{Two empirical data sets satisfying relation given in \eqref{eq:general}. Panel a) The empirical time series $X_1^1, X_2^1, \ldots X_{20000}^1$.  Panel b) The empirical time series $X_1^2, X_2^2, \ldots X_{40000}^2$  with the corresponding estimates of the critical changes points $\hat{l}_{V}$ (red dashed line) and $\hat{l}_{C}$ (olive dashed line) obtained by using approaches described in Sections ~\ref{sec:scaleBased} and  ~\ref{sec:GajdaEtAll}}
\label{fig:real_timeseries_stats}
\end{figure}

In this section we apply the discussed techniques to the real data set presented in Fig.~\ref{fig:real_timeseries}. The data represent the vibration signal from a heavy duty crusher used in a mineral processing plant. The
crusher is a kind of machine which use a metal surface to crumble materials into small fractional
pieces. During this process, as well as during entering material stream into the crusher, a lot of
impacts/shocks appear. They are present in vibration signal acquired from the bearings housing \cite{vibro}. In the signal we observe that some statistical properties changes over time.  The empirical analysis is built on top of two different data sets manually extracted from the original time series, see Fig.~\ref{fig:real_timeseries}  panel b) and panel c). We denote them $X_1^1, X_2^1, \ldots X_{20000}^1$ and $X_1^2, X_2^2, \ldots X_{40000}^2$, respectively. For both cases, one may easily observe that change in the scale parameter of the driving distribution has taken place. Therefore we can suspect that the property given in \eqref{eq:general} holds for each of the observed time series. At the first stage of our analysis, we discuss the results of estimation procedures for critical change points ($\hat{l}_V$ and $\hat{l}_C$) calculated based on $V$ and $C$ statistics formulated in \eqref{eq:NewStats} and \eqref{eq:GajdaEtAllStats}, respectively. For ease of reading, the results are covered in Table~\ref{tbl:real_data_resuts}. Moreover the obtained estimates of the critical points $\hat{l}_C$ and $\hat{l}_V$ along with reference statistics $C$ and $V$ are depicted in Fig.~\ref{fig:real_timeseries_stats}.

\noindent For the time series consisting of $20000$ observations the estimates of the critical change points $\hat{l}_V$ and $\hat{l}_C$ are $6282$ and $7591$.
The corresponding $p$-values of the performed statistical tests given in \eqref{eq:ansariStats} and \eqref{eq:quantileTestStats} are  $1.15\cdot 10^{-25}$ and $1$ allowing us to reach conclusions, where the null hypothesis $\mathcal{H}_0$ is rejected in case of the estimate provided by $V$, whereas there is not enough evidence to reject the null hypothesis $\mathcal{H}_0$ in case of statistics $C$ at the $0.05$ significance level. The similar conclusions are made for the second time series that consists of $40000$ observations as the obtained $p$-values for $\hat{l}_V = 14389$ and $\hat{l}_C=12565$ estimates are $3.445\cdot 10 ^{-15}$ and $1$, respectively.
\begin{table}[h!]
\centering
\begin{tabular}{ccc}
\centering
& Sample I & Sample II \\
Observations & $20\ 000$ & $40\ 000$ \\\hline\hline
$\hat{l}_C$ & $7591\ (\approx 1)$ & $12565\ (\approx 1)$\\\hline
$\hat{l}_V$ & $6282\ (<0.01)$ & $14389\ (<0.01)$\\\hline
\end{tabular}
\caption{Results of the performed analysis for estimating critical change points $\hat{l}_C$ and $\hat{l}_V$. The corresponding $p$-values are given in parentheses.}
\label{tbl:real_data_resuts}
\end{table}

\begin{figure}[h!]
\centering
\includegraphics[scale=0.35]{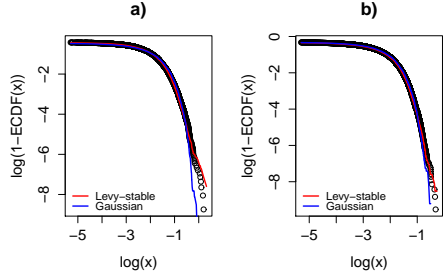}
\caption{The right tail fit to the empirical cumulative distribution function for the corresponding sub samples extracted from $X_1^1, X_2^1, \ldots,X_{20000}^1$ given on a double logarithmic scale.
Panel a) Empirical sample $X_1^1,\ldots, X_{6282}^1$. Panel b) Empirical sample  $X_{6283}^1,\ldots, X_{20000}^1$. Both panels include L\'evy-Stable (red) and Gaussian (blue) fits to the corresponding empirical cumulative distribution functions produced by the real data (black circles).}
\label{fig:real_data_tails_20000}
\end{figure}

\begin{figure}[h!]
\centering
\includegraphics[scale=0.35]{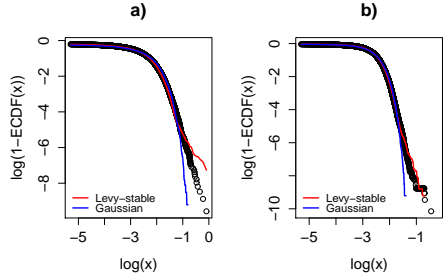}
\caption{The right tail fit to the empirical cumulative distribution function for the corresponding sub samples extracted from $X_1^2, X_2^2, \ldots,X_{40000}^2$ given on a double logarithmic scale.
Panel a) Empirical sample $X_1^2,\ldots, X_{14389}^2$. Panel b) Empirical sample  $X_{14390}^1,\ldots, X_{40000}^2$. Both panels include L\'evy-Stable (red) and Gaussian (blue) fits to the corresponding empirical cumulative distribution functions produced by the real data (black circles).}
\label{fig:real_data_tails_40000}
\end{figure}

\noindent The complementary part of this section is to provide details about potential family of distributions for the sub samples obtained by splitting the original sets of observations $X_1^1, X_2^1, \ldots X_{20000}^1$ and $X_1^2, X_2^2, \ldots X_{40000}^2$ into two parts according to the rule provided by the appropriate estimates of the critical change points $\hat{l}_V$ given in Table~\ref{tbl:real_data_resuts}. As a consequence the first time series $X_1^1, X_2^1, \ldots X_{20000}^1$ has been divided into $X_1^1, X_2^1, \ldots X_{6282}^1$ and $X_{6283}^1, X_{6284}^1, \ldots X_{20000}^1$, accordingly. Furthermore two different samples consisting of $X_{1}^2, X_{2}^2, \ldots X_{14389}^2$ and $X_{14390}^2, X_{14391}^2, \ldots X_{40000}^2$ have been created based on the second time series $X_1^2, X_2^2, \ldots X_{40000}^2$.
For the purposes of finding appropriate distribution that fits well to the empirical data, we follow the procedure where L\'evy-stable distribution $\mathcal{S}(\alpha,\beta,\sigma,\mu)$ and Gaussian distribution $\mathcal{N}(\mu,\sigma)$ are considered. In order to fit the parameters of L\'evy-stable distribution $\mathcal{S}(\alpha,\beta,\sigma,\mu)$ we apply method described in \cite{Koutrouvelis1980}. The results of the estimated parameters of
L\'evy-stable distribution $\mathcal{S}(\alpha,\beta,\sigma,\mu)$ and Gaussian distribution $\mathcal{N}(\mu,\sigma)$ for all samples are presented in Table~\ref{tbl:estSample1_1}. Moreover L\'evy-stable and Gaussian right tail fits on a double logarithmic scale for samples $X_1^1, X_2^1, \ldots X_{20000}^1$ and $X_1^2, X_2^2, \ldots X_{40000}^2$ are presented in Fig.~\ref{fig:real_data_tails_20000} and~\ref{fig:real_data_tails_40000}, respectively.
\begin{table}[h!]
\centering
\begin{tabular}{cccccccc}
Sample&Distribution fit & $\alpha$ & $\beta$ & $\mu$ & $\sigma$ & KS & AD  \\\hline\hline
\multirow{2}{*}{$X_1^1,\ldots,X_{6282}^1$}&$\mathcal{S}(\alpha,\beta,\sigma,\mu)$ & $1.8587$ & $-0.1110$ & $0.0769$ & $0.1370$ & $ 0.0097\ (0.9274)$ & $0.6889\ (0.5097)$\\\cline{2-8}
&$\mathcal{N}(\mu,\sigma)$ & & & $0.0789$ & $0.2165$ & $0.3448\ (<0.05)$ & $13.0953\ (<0.05)$\\\hline
\multirow{2}{*}{$X_{6283}^1,\ldots,X_{20000}^1$}&$\mathcal{S}(\alpha,\beta,\sigma,\mu)$ & $1.9429$ & $-0.1160$ & $0.0769$ & $0.0887$ & $ 0.0069\ (0.8963)$ & $0.6889\ (0.5097)$\\\cline{2-8}
&$\mathcal{N}(\mu,\sigma)$ & & & $0.0771$ & $0.1314$ & $0.4023\ (<0.05)$ & $7.7176\ (<0.05)$\\\hline\hline
\multirow{2}{*}{$X_{1}^2,\ldots, X_{14389}^2$}& $\mathcal{S}(\alpha,\beta,\sigma,\mu)$ & $1.8190$ & $-0.0592$ & $0.0577$ & $0.0763$ & $  0.0085\ (0.6774)$ & $0.6889\ (0.5097)$\\\cline{2-8}
&$\mathcal{N}(\mu,\sigma)$ & & & $0.0766$ & $0.0951$ & $ 0.4296\ (<0.05)$ & $\infty\ (<0.05)$\\\hline
\multirow{2}{*}{$X_{14390}^2,\ldots,X_{40000}^2$}&$\mathcal{S}(\alpha,\beta,\sigma,\mu)$ & $1.8814 $ & $0.0610$ & $0.0281$ & $0.0797$ & $ 0.0064\ (0.6761)$ & $0.6889\ (0.5097)$\\\cline{2-8}
& $\mathcal{N}(\mu,\sigma)$ & & & $0.0797$ & $0.0443$ & $0.4801\ (<0.05)$ & $\infty\ (<0.05)$\\\hline\hline
\end{tabular}
\caption{Estimated parameters of L\'evy-stable distribution $\mathcal{S}(\alpha,\beta,\sigma,\mu)$ and Gaussian distribution $\mathcal{N}(\mu,\sigma)$ for $X_1^1, X_2^1, \ldots,X_{20000}^1$ and $X_1^2, X_2^2, \ldots,X_{40000}^2$ samples as well as the KS and AD statistics for L\'evy stable distribution and corresponding p-values (presented in brackets).}
\label{tbl:estSample1_1}
\end{table}

\noindent In terms of the fits to the first empirical data sets  $X_1^1,\ldots, X_{6282}^1$ and $X_{6283}^1,\ldots, X_{20000}^1$, we observe that the obtained estimates of L\'evy-stable distribution $\mathcal{S}(\alpha,\beta,\sigma,\mu)$ provide better fits to the data than the corresponding Gaussian distribution $\mathcal{N}(\mu,\sigma)$, see Fig.~\ref{fig:real_data_tails_20000}. The observed superiority of the L\'evy-stable distribution is pronounced by the fact, that in both cases the calculated Anderson-Darling and Kolmogorov-Smirnov test statistics \cite{bw,tests1}, force us to reject the hypothesis that the data is normally distributed in favour of L\'evy-stable distribution.

\noindent To be more precise, for the first vector of observations $X_1^1,\ldots, X_{6282}^1$ the values of the Anderson-Darling and Kolmogorov statistics produced by L\'evy-stable distribution are $0.6889$ and $0.0097$, respectively. The corresponding $p$-values are $0.9274$ and $0.5097$ allowing us to accept the L\'evy-stable law a model of the considered set of data. The values of the test statistics for the Gaussian fit model yield $p$-values of less than $0.05$ forcing us to reject the assumption that the data is normally distributed at the $0.05$ significance level. The similar conclusions are made for the second time series $X_{6283}^1,\ldots, X_{20000}^1$ as the obtained values
of the tests statistics satisfy the same relationship as in case of the aforementioned sample $X_1^1,\ldots, X_{6282}^1$, see Table~\ref{tbl:estSample1_1}. Subsequently, we apply the same technique to the second set of empirical observations, i.e. $X_{1}^2, X_{2}^2, \ldots X_{14389}^2$ and $X_{14390}^2, X_{14391}^2, \ldots X_{40000}^2$. It appears that L\'evy-stable distributions fit very well to the both samples, yielding relative small values of the Anderson-Darling ($0.6889$ and $0.6889$) and Kolmogorow-Smirnov ($0.0085$ and $0.0064$) test statistics for the samples, $X_{1}^2, X_{2}^2, \ldots X_{14389}^2$ and $X_{14390}^2, X_{14391}^2, \ldots X_{40000}^2$, respectively.
The associated $p$-values for the first sample  $X_{1}^2, X_{2}^2, \ldots X_{14389}^2$ are $0.5097$ and $0.6774$, whereas for the second sample $X_{14390}^2, X_{14391}^2, \ldots X_{40000}^2$ are equal to $0.5097$ and $0.6761$, respectively, giving us an opportunity to accept the hypothesis that both samples are distributed according to the appropriate L\'evy-stable law.

\noindent The results obtained by fitting the Gaussian model to the both samples force us to reject the hypothesis that the data is normally distributed as the observed Anderson-Darling and  Kolmogorow-Smirnow test produces $p$-vales of less than $0.05$ for all considered cases. The corresponding tail fits of the empirical cumulative distributive function empirical for samples $X_{1}^2, X_{2}^2, \ldots X_{14389}^2$ and $X_{14390}^2, X_{14391}^2, \ldots X_{40000}^2$ are presented in Fig.~\ref{fig:real_data_tails_40000}.
\section{Conclusions}
In this paper we have introduced a new technique of time series  segmentation, i.e extraction from the originally signal such parts that have similar statistical properties.   The proposed technique is based on the simple statistics (absolute deviation about the median) based on given time series. The structural break point detection method proposed here uses the  Adaptive Regression Splines technique, which is an extension of the classical regression.  Except the introduction of the   estimation procedure for the recognition of the critical point that divides the observed time series into two
regimes, we have also developed the statistical test of testing  two-regimes behavior. The
universality of the presented methodology comes from the fact that it does
not assume the distribution of the examined time series, therefore, it can
be applied to a rich class of real data sets.   The theoretical results we have illustrated using the simulated
time series and analysis of  real data set  that represents the vibration signal from a heavy duty crusher used in a mineral processing plant.

\end{document}